\documentclass{article}
\usepackage[final]{graphicx}   
\usepackage{psfrag}   
\usepackage{amsmath,amsfonts,amssymb,amsxtra}

\newcommand {\Real}{\ensuremath{{\mathbb{R}}}}

\newcommand {\Complex}{\ensuremath{{\mathbb{C}}}}

\newcommand{\A}{\ensuremath{\mathcal A}}

\newcommand{\N}{\ensuremath{\mathcal N}}

\newcommand{\vi}{\ensuremath{{\mathbf{v}}}}
\newcommand{\ex}{\ensuremath{{\mathbf{x}}}}

\newcommand{\one}{\ensuremath{{\mathbf{1}}}}

\newtheorem{theorem}{Theorem}

\newtheorem{lemma}{Lemma}

\newtheorem{fact}{Fact}

\newenvironment{proof}{\noindent {\bf Proof.}}{\hfill \hspace*{1pt}\hfill$\blacksquare$}

\begin{document}
\title{LQR-based coupling gain for synchronization of linear systems} 
\author{S. Emre Tuna\\
{\small {\tt tuna@eee.metu.edu.tr}} }
\maketitle

\begin{abstract}
Synchronization control of coupled continuous-time linear systems is
studied.  For identical systems that are stabilizable, a linear
feedback law obtained via algebraic Riccati equation is shown to
synchronize any fixed directed network of any number of coupled
systems provided that the coupling is strong enough. The strength of
coupling is determined by the smallest distance of a nonzero
eigenvalue of the coupling matrix to the imaginary axis. A dual
problem where detectable systems that are coupled via their outputs is
also considered and solved.
\end{abstract}

\section{Introduction}
Synchronization is, at one hand, a desired behaviour in many dynamical
systems related to numerous technological applications
\cite{fabiny93,fax04,cortes04}; and, at the other, a
frequently-encountered phenomenon in biology
\cite{stoop00,walker69,aldridge76}. On top of that, it is an
important system theoretical topic on its own right.  For instance,
there is nothing keeping us from seeing a simple Luenberger observer
\cite{luenberger64} with decaying error dynamics as a two-agent system
where the agents globally synchronize. Since people lack anything but
good reasons to investigate synchronization, a wealth of literature
has been formed, mostly in recent times
\cite{boccaletti06,strogatz01,wang02}.

An essential problem from a control theory point of view is to find
conditions that imply synchronization of a number of coupled
individual systems. This problem, which is usually studied under the
name {\em synchronization stability}, has been attacked by many and
from various angles. Two cases are of particular interest: (i) where
the dynamics of individual systems are primitive (such as that of an
integrator) yet the coupling between them is considered time-varying;
and (ii) where the individual systems are let be more sophisticated,
however, coupled via a fixed interconnection. The studies concentrated
on the first case have resulted in the emergence of the area now known
as {\em consensus} in {\em multi-agent systems}
\cite{olfati07,moreau05,cortes07,jadbabaie03,
ren05,angeli06,blondel05,tsitsiklis86} where fairly weak conditions on
the interconnection have been established under which the states of
individual systems converge to a common point that is fixed in
space. The second case, into which the problem studied in this note
falls, has also accommodated important theoretical developments
especially by using tools from algebraic graph theory
\cite{wu05con}. Using Lyapunov functions, it has been shown that
spectrum of the coupling matrix plays a crucial role in determining
the stability of synchronization
\cite{wu95,pecora98} notwithstanding it need not necessarily be
explicitly known \cite{belykh06}. It has also been shown that
passivity theory can be useful in studying stability
provided that the interconnection is symmetric
\cite{arcak07,pogromsky01,stan07}.

In this note we consider identical individual system dynamics that
are linear time-invariant $\dot{x}=Ax+Bu$. Under the weakest possible
assumption that pair $(A,\,B)$ is stabilizable, we search for a
feedback law $\kappa(A,\,B)$ which would guarantee asymptotic
synchronization for any fixed (directed) interconnection of arbitrary
number of coupled systems provided that the coupling is strong
enough. Following the tradition, we use the spectral information of
the coupling matrix to determine the strength of connectedness.
Namely, the farther the second eigenvalue (with largest real part)
from the imaginary axis the more connected the network. We show that a
linear $\kappa$ solving a linear quadratic regulation (LQR) problem
performs the task.

It is worth noting that most of the existing work on synchronization
focuses on analysis rather than design. Also, the individual system
dynamics are usually taken to be stable (i.e. the trajectories of the
uncoupled system are required to be bounded.) In those respects the
issue we deal with in this work is relatively different. Of particular
relevance to this note are the works \cite{aut7874} and
\cite{wu05}. The former provides, for linear individual system
dynamics, a linear feedback law that guarantees synchronization for
all connected interconnections (regardless of the strength of
coupling) under the extra assumption that matrix $A$ is neutrally
stable. The latter establishes sufficient conditions for nonlinear
individual system dynamics so that synchronization is achieved for
strong enough coupling.

In the remainder of the paper we first provide notation and
definitions.  Then, in Section~\ref{sec:problem}, we formalize the
problem and state our objectives. In Section~\ref{sec:main} we show
via our main theorem that optimal control theory yields us a feedback
law which serves our purpose, i.e. synchronizes coupled systems for
all network topologies with strong enough coupling. Finally, in
Section~\ref{sec:dual}, a dual problem is formulated and solved.

\section{Notation and definitions}

Let $\Real_{\geq 0}$ denote set of nonnegative real numbers and
$|\cdot|$ 2-norm. For $\lambda\in\Complex$ let ${\rm Re}(\lambda)$
denote the real part of $\lambda$. Identity matrix in $\Real^{n\times
n}$ is denoted by $I_{n}$ and zero matrix in $\Real^{m\times n}$ by
$0_{m\times n}$. Conjugate transpose of a matrix $A$ is denoted by
$A^{H}$. Matrix $A\in\Complex^{n\times n}$ is {\em Hurwitz} if all
of its eigenvalues have strictly negative real parts.\footnote{Note
that $A$ is Hurwitz if there exists a symmetric positive definite
matrix $P$ such that $A^{H}P+PA<0$.} Given $B\in\Real^{n\times m}$,
$C\in\Real^{m\times n}$, and $A\in\Real^{n\times n}$; pair $(A,\,B)$
is {\em stabilizable} if there exists $K\in\Real^{m\times n}$ such
that $A-BK$ is Hurwitz; pair $(C,\,A)$ is {\em detectable} if
$(A^{T},\,C^{T})$ is stabilizable. Let $\one\in\Real^{p}$ denote the
vector with all entries equal to one.

{\em Kronecker product} of $A\in\Complex^{m\times n}$ and $B\in\Complex^{p\times q}$ is
\begin{eqnarray*}
A\otimes B:=
\left[
\begin{array}{ccc}
a_{11}B & \cdots & a_{1n}B\\
\vdots  & \ddots & \vdots\\
a_{m1}B & \cdots & a_{mn}B
\end{array}
\right]
\end{eqnarray*}
Kronecker product comes with the property $(A\otimes
B)(C\otimes D)=(AC)\otimes(BD)$ (provided that products $AC$ and $BD$
are allowed.) 

A ({\em directed}) {\em graph} is a pair $(\N,\,\A)$ where $\N$ is a
nonempty finite set (of {\em nodes}) and $\A$ is a finite collection
of pairs ({\em arcs}) $(n_{i},\,n_{j})$ with $n_{i},\,n_{j}\in\N$. A
{\em path} from $n_{1}$ to $n_{\ell}$ is a sequence of nodes
$\{n_{1},\,n_{2},\,\ldots,\,n_{\ell}\}$ such that $(n_{i},\,n_{i+1})$
is an arc for $i\in\{1,\,2,\,\ldots,\,\ell-1\}$. A graph is {\em
connected} if it has a node to which there exists a path from every
other node.\footnote{Note that this definition of connectedness for
directed graphs is weaker than strong connectivity and stronger than
weak connectivity. In \cite{wu05} an equivalent condition (for
connectedness) is given as that the graph contains a spanning tree.}

The graph of a matrix $\Gamma:=[\gamma_{ij}]\in\Real^{p\times p}$ is
the pair $(\N,\,\A)$ where $\N =\{n_{1},\,n_{2},\,\ldots,\,n_{p}\}$
and $(n_{i},\,n_{j})\in\A$ iff $\gamma_{ij}>0$. Matrix $\Gamma$ is
said to be {\em connected} if it
satisfies:
\begin{enumerate}
\item[(i)] $\gamma_{ij}\geq 0$ for $i\neq j$;
\item[(ii)] each row sum equals 0;
\item[(iii)] its graph is connected.
\end{enumerate}

A connected $\Gamma$ has an eigenvalue at $\lambda=0$ with eigenvector
$\one$, i.e. $\Gamma\one=0$, and all its other eigenvalues have real
parts strictly negative.\footnote{For the sake of completeness we
provide a proof of this fact in Appendix. }
When we write ${\rm Re}(\lambda_{2}(\Gamma))$ we mean the real part of
a nonzero eigenvalue of $\Gamma$ closest to the imaginary axis.

Given maps $\xi_{i}:\Real_{\geq 0}\to\Real^{n}$ for
$i\in\{1,\,2,\,\ldots,\,p\}$ and a map $\bar\xi:\Real_{\geq
0}\to\Real^{n}$, the elements of the set
$\{\xi_{i}(\cdot):i=1,\,2,\,\ldots,\,p\}$ are said to {\em synchronize
to} $\bar{\xi}(\cdot)$ if $|\xi_{i}(t)-\bar\xi(t)|\to 0$ as
$t\to\infty$ for all $i$. The elements of the set
$\{\xi_{i}(\cdot):i=1,\,2,\,\ldots,\,p\}$ are said to synchronize if
they synchronize to some $\bar{\xi}(\cdot)$.

\section{Problem statement}\label{sec:problem}

\subsection{Systems under study}
We consider $p$ identical linear systems
\begin{eqnarray}\label{eqn:system}
{\dot{x}}_{i}=Ax_{i}+Bu_{i}\, ,\quad i=1,\,2,\,\ldots,\,p
\end{eqnarray}
where $x_{i}\in\Real^{n}$ is the {\em state} and $u_{i}\in\Real^{m}$ is the {\em input}
of the $i$th system. Matrices
$A$ and $B$ are of proper dimensions. The solution of $i$th system at
time $t\geq 0$ is denoted by $x_{i}(t)$. In this paper we consider
the case where at each time instant (only) the following information
\begin{eqnarray}\label{eqn:z}
z_{i}&=&\sum_{j=1}^{p}\gamma_{ij}(x_{j}-x_{i})
\end{eqnarray}
is available to $i$th system to determine an input value where
$\gamma_{ij}$ are the entries of the matrix $\Gamma\in\Real^{p\times p}$ describing
the network topology. Nondiagonal entries of $\Gamma$ are nonnegative
and each row sums up to zero. That is, the coupling between systems is
diffusive.

\subsection{Assumptions made}
We only make the following assumption on systems~\eqref{eqn:system}
which will henceforth hold.

{\bf (A1)} Pair $(A,\,B)$ is stabilizable.

\subsection{Objectives}
We have two objectives in this paper.

{\bf (O1)} Show that for each $\delta>0$ there exists a linear
feedback law $K\in\Real^{m\times n}$ such that, for all $p$ and
connected $\Gamma\in\Real^{p\times p}$ with $-{\rm
Re}(\lambda_{2}(\Gamma))\geq\delta$, solutions of
systems~\eqref{eqn:system} for $u_{i}=Kz_{i}$, where $z_{i}$ is as in
\eqref{eqn:z}, globally (i.e. for all initial conditions)
synchronize. 

{\bf (O2)} Compute one such $K$.

\section{Main result}\label{sec:main}

For later use in this section we first borrow a well-known result from
optimal control theory \cite{sontag98}: Given a stabilizable pair
$(A,\,B)$, where $A\in\Real^{n\times n}$ and $B\in\Real^{n\times m}$,
the following algebraic Riccati equation
\begin{eqnarray}\label{eqn:riccati}
A^{T}P+PA+I_{n}-PBB^{T}P=0
\end{eqnarray}  
has a (unique) solution $P=P^{T}>0$. One can rewrite
\eqref{eqn:riccati} as
\begin{eqnarray*}
(A-BB^{T}P)^{T}P+P(A-BB^{T}P)+(I_{n}+PBB^{T}P)=0
\end{eqnarray*}
whence we infer that $A-BB^{T}P$ is Hurwitz. 

\begin{lemma}\label{lem:lqr}
Let $A\in\Real^{n\times n}$ and $B\in\Real^{n\times m}$ satisfy
\eqref{eqn:riccati} for some symmetric positive definite $P$. Then for
all $\sigma\geq 1$ and $\omega\in\Real$ matrix
$A-(\sigma+j\omega)BB^{T}P$ is Hurwitz.
\end{lemma}

\begin{proof}
Let $\varepsilon:=\sigma-1\geq 0$. Write 
\begin{eqnarray}\label{eqn:complexlyap}
\lefteqn{(A-(\sigma+j\omega)BB^{T}P)^{H}P+P(A-(\sigma+j\omega)BB^{T}P)}\nonumber\\
&&=(A-(\sigma-j\omega)BB^{T}P)^{T}P+P(A-(\sigma+j\omega)BB^{T}P)\nonumber\\
&&=(A-(1+\varepsilon)BB^{T}P)^{T}P+P(A-(1+\varepsilon)BB^{T}P)\nonumber\\
&&=(A-BB^{T}P)^{T}P+P(A-BB^{T}P)-2\varepsilon PBB^{T}P\nonumber\\
&&=-I_{n}-(1+2\varepsilon)PBB^{T}P\,.
\end{eqnarray}
Finally, observe that
\eqref{eqn:complexlyap} is nothing but (complex) Lyapunov equation.
\end{proof}

Below is our main result.
\begin{theorem}\label{thm:main}
Consider systems \eqref{eqn:system}. Let $K:=B^{T}P$ where $P$ is the
solution to \eqref{eqn:riccati}. Given $\delta>0$, for all $p$ and
connected $\Gamma\in\Real^{p\times p}$ with $-{\rm
Re}(\lambda_{2}(\Gamma))\geq\delta$, solutions $x_{i}(\cdot)$ for
$i=1,\,2,\,\ldots,\,p$ and $u_{i}=\max\{1,\,\delta^{-1}\}Kz_{i}$,
where $z_{i}$ is as in \eqref{eqn:z}, globally synchronize to
\begin{eqnarray*}
\bar{x}(t):=
(r^{T}\otimes e^{At})
\left[
\begin{array}{c}
x_{1}(0)\\
\vdots\\
x_{p}(0)
\end{array}
\right]
\end{eqnarray*}
where $r\in\Real^{p}$ satisfies $r^{T}\Gamma=0$ and $r^{T}\one=1$.
\end{theorem}

\begin{proof}
Let $K_{\delta}:=\max\{1,\,\delta^{-1}\}K$. Combine \eqref{eqn:system}
and \eqref{eqn:z} to obtain
\begin{eqnarray}\label{eqn:couplexi}
\dot{x}_{i}=Ax_{i}+BK_{\delta}\sum_{j=1}^{p}\gamma_{ij}(x_{j}-x_{i})\,.
\end{eqnarray}
Stack individual system states as
$\ex:=[x_{1}^{T}\ x_{2}^{T}\ \ldots\ x_{p}^{T}]^{T}$. Then we can express 
\eqref{eqn:couplexi} as
\begin{eqnarray}\label{eqn:couplex}
\dot{\ex}=(I_{p}\otimes A + \Gamma\otimes BK_{\delta})\ex\,.
\end{eqnarray}
Now let $Y\in\Complex^{p\times (p-1)}$, $W\in\Complex^{(p-1)\times p}$,
$V\in\Complex^{p\times p}$, and upper triangular
$\Delta\in\Complex^{(p-1)\times(p-1)}$ be such that
\begin{eqnarray*}
V=\left[\one\ Y\right]\ ,\quad 
V^{-1}=
\left[
\begin{array}{c}
r^{T}\\
W
\end{array}
\right]
\end{eqnarray*} 
and 
\begin{eqnarray*}
V^{-1}\Gamma V=
\left[
\begin{array}{cc}
0&0\ \cdots\ 0\\
\begin{array}{c}
0\\
\vdots\\
0
\end{array}
& \Delta
\end{array}
\right]
\end{eqnarray*}
Note that the diagonal entries of $\Delta$ are nothing but the nonzero
eigenvalues of $\Gamma$ which we know have real parts no greater than
$-\delta$. Engage the change of variables $\vi:=(V^{-1}\otimes
I_{n})\ex$ and modify
\eqref{eqn:couplex} first into
\begin{eqnarray*}
\dot{\vi}=(I_{p}\otimes A+V^{-1}\Gamma V\otimes BK_{\delta})\vi
\end{eqnarray*}
and then into
\begin{eqnarray}\label{eqn:vi}
\dot{\vi}=
\left[
\begin{array}{cc}
A&0_{n\times(p-1)n}\\
0_{(p-1)n\times n}& I_{p-1}\otimes A+\Delta\otimes BK_{\delta}
\end{array}
\right]\vi
\end{eqnarray}
Observe that $I_{p-1}\otimes A+\Delta\otimes BK_{\delta}$ is upper block
triangular with (block) diagonal entries of the form $A+\lambda_{i}
BK_{\delta}$ for $i=2,\,3,\,\ldots,\,p$ with ${\rm Re}(\lambda_{i})\leq
-\delta$. Lemma~\ref{lem:lqr} implies therefore that $I_{p-1}\otimes
A+\Delta\otimes BK_{\delta}$ is Hurwitz. Thus \eqref{eqn:vi} implies
\begin{eqnarray*}
\left|\vi(t)-\left[
\begin{array}{cc}
e^{At}&0_{n\times(p-1)n}\\
0_{(p-1)n\times n}& 0_{(p-1)n\times(p-1)n}
\end{array}
\right]\vi(0)\right|\to 0
\end{eqnarray*} 
as $t\to\infty$ which yields
\begin{eqnarray*}
\left|\ex(t)-(\one r^{T}\otimes e^{At})\ex(0)\right|\to 0\,.
\end{eqnarray*} 
Hence the result.
\end{proof}

By Theorem~\ref{thm:main} we attain our objectives. In the next
section we provide a dual result which may be more useful in certain
applications.

\section{Dual problem}\label{sec:dual}
Let $p$ identical linear systems be 
\begin{eqnarray}\label{eqn:systemdual}
\dot{x}_{i}=A^{T}x_{i}+u_{i}\, ,\quad y_{i}=B^{T}x_{i}\, ,\quad i=1,\,2,\,\ldots,\,p
\end{eqnarray}
where $x_{i}\in\Real^{n}$ is the state, $u_{i}\in\Real^{n}$ is the
input, and $y_{i}\in\Real^{m}$ is the {\em output} of the $i$th
system. Matrices $A^{T}$ and $B^{T}$ are of proper dimensions and they
make a detectable pair $(B^{T},\,A^{T})$. Now consider the case where
at each time instant the following information
\begin{eqnarray}\label{eqn:zdual}
z_{i}&=&\sum_{j=1}^{p}\gamma_{ij}(y_{j}-y_{i})
\end{eqnarray}
is available to $i$th system to determine an input value. Not
surprisingly, the following result accrues.

\begin{theorem}
Consider systems \eqref{eqn:systemdual}. Let $L:=PB$ where $P$ is the
solution to \eqref{eqn:riccati}. Given $\delta>0$, for all $p$ and
connected $\Gamma\in\Real^{p\times p}$ with $-{\rm
Re}(\lambda_{2}(\Gamma))\geq\delta$, solutions $x_{i}(\cdot)$ for
$i=1,\,2,\,\ldots,\,p$ and $u_{i}=\max\{1,\,\delta^{-1}\}Lz_{i}$,
where $z_{i}$ is as in \eqref{eqn:zdual}, globally synchronize to
\begin{eqnarray*}
\bar{x}(t):=
(r^{T}\otimes e^{A^{T}t})
\left[
\begin{array}{c}
x_{1}(0)\\
\vdots\\
x_{p}(0)
\end{array}
\right]
\end{eqnarray*}
where $r\in\Real^{p}$ satisfies $r^{T}\Gamma=0$ and $r^{T}\one=1$.
\end{theorem}

\section{Conclusion}
For identical (unstable) linear systems, we have shown that no more
than stabilizability (detectability) is necessary for an LQR-based
linear feedback law to exist under which coupled systems synchronize
for all network topologies with strong enough coupling. We have
considered directed networks and measured the strength of coupling via
the real part of a nonzero eigenvalue (of the coupling matrix) closest
to the imaginary axis.

\appendix
\section{Proof of a fact}

\begin{fact}
A connected $\Gamma\in\Real^{p\times p}$ has an eigenvalue at
$\lambda=0$ with eigenvector $\one$ and all
its other eigenvalues have strictly negative real parts.
\end{fact}

\begin{proof}
First part of the statement directly follows from the
definition. Since the entries of each row of $\Gamma$ sum up to zero,
it trivially follows that $\Gamma\one=0$.

Consider the dynamical system
\begin{eqnarray*}
\dot{x}=\Gamma{x}
\end{eqnarray*}
where $x\in\Real^{p}$ with entries $x_{i}\in\Real$ for $i=1,\,2,\,\ldots,\,p$. 
Since $\Gamma\one=0$ we can write
\begin{eqnarray*}
\dot{x}_{i}=\sum_{j=1}^{p}\gamma_{ij}(x_{j}-x_{i})\,.
\end{eqnarray*} 
By definition $\gamma_{ij}\geq 0$ for $i\neq j$. That brings us the
following.  For each $\tau\geq 0$,
$x_{i}(t)\in[\min_{i}x_{i}(\tau),\,\max_{i}x_{i}(\tau)]$ for all $i$
and $t\geq\tau$. In addition, since the graph of $\Gamma$ is
connected, interval $[\min_{i}x_{i}(t),\,\max_{i}x_{i}(t)]$ must
shrink to a point as $t\to\infty$. In other words, there exists
$\bar{x}\in[\min_{i}x_{i}(0),\,\max_{i}x_{i}(0)]$ such that
\begin{eqnarray}\label{eqn:app2}
\lim_{t\to\infty} x_{i}(t)=\bar{x}
\end{eqnarray}
for all $i$ (see \cite{moreau04}.) That readily implies by simple
stability arguments that if $\lambda\in\Complex$ is an eigenvalue of
$\Gamma$ then ${\rm Re}(\lambda)\leq 0$. Also, for ${\rm
Re}(\lambda)=0$, the size of the associated Jordan block cannot be
greater than unity. Possibility of a purely imaginary eigenvalue means
sustaining oscillations and is therefore ruled out by
\eqref{eqn:app2}. The only case left unconsidered is a
second eigenvalue at the origin. That would imply, since it cannot be
related to a Jordan block of size two or greater, that there exists a
nonzero (eigen)vector $v\in\Real^{p}$ such that $\Gamma{v}=0$ and
$v\neq\alpha\one$ for any $\alpha\in\Real$. That contradicts
\eqref{eqn:app2}, too.
\end{proof} 

\bibliographystyle{plain}
\bibliography{references} 

\begin{thebibliography}{10}

\bibitem{aldridge76}
J.~Aldridge and E.K. Pye.
\newblock Cell density dependence of oscillatory metabolism.
\newblock {\em Nature}, 259:670--671, 1976.

\bibitem{angeli06}
D.~Angeli and P.-A. Bliman.
\newblock Stability of leaderless discrete-time multi-agent systems.
\newblock {\em Mathematics of Control, Signals \& Systems}, 18:293--322, 2006.

\bibitem{arcak07}
M.~Arcak.
\newblock Passivity as a design tool for group coordination.
\newblock {\em IEEE Transactions on Automatic Control}, 52:1380--1390, 2007.

\bibitem{belykh06}
I.~Belykh, V.~Belykh, and M.~Hasler.
\newblock Generalized connection graph method for synchronization in
  asymmetrical networks.
\newblock {\em Physica D}, 224:42--51, 2006.

\bibitem{blondel05}
V.D. Blondel, J.M. Hendrickx, A.~Olshevsky, and J.N. Tsitsiklis.
\newblock Convergence in multiagent coordination, consensus, and flocking.
\newblock In {\em Proc. of the 44th IEEE Conference on Decision and Control},
  pages 2996--3000, 2005.

\bibitem{boccaletti06}
S.~Boccaletti, V.~Latora, Y.~Moreno, M.~Chavez, and D.-U. Hwang.
\newblock Complex networks: structure and dynamics.
\newblock {\em Physics Reports-Review Section of Physics Letters},
  424:175--308, 2006.

\bibitem{cortes07}
J.~Cortes.
\newblock Distributed algorithms for reaching consensus on general functions.
\newblock {\em Automatica}, 2007.

\bibitem{cortes04}
J.~Cortes, S.~Martinez, T.~Karatas, and F.~Bullo.
\newblock Coverage control for mobile sensing networks.
\newblock {\em IEEE Transactions on Robotics and Automation}, 20:243--255,
  2004.

\bibitem{fabiny93}
L.~Fabiny, P.~Colet, R.~Roy, and D.~Lenstra.
\newblock Coherence and phase dynamics of spatially coupled solid-state lasers.
\newblock {\em Physical Review A}, 47:4287--4296, 1993.

\bibitem{fax04}
J.A. Fax and R.M. Murray.
\newblock Information flow and cooperative control of vehicle formations.
\newblock {\em IEEE Transactions on Automatic Control}, 49:1465--1476, 2004.

\bibitem{jadbabaie03}
A.~Jadbabaie, J.~Lin, and A.S. Morse.
\newblock Coordination of groups of mobile autonomous agents using nearest
  neighbor rules.
\newblock {\em IEEE Transactions on Automatic Control}, 48:988--1001, 2003.

\bibitem{luenberger64}
D.G. Luenberger.
\newblock Observing the state of a linear system.
\newblock {\em IEEE Transactions on Military Electronics}, pages 74--80, April
  1964.

\bibitem{moreau04}
L.~Moreau.
\newblock Stability of continuous-time distributed consensus algorithms.
\newblock {\em \tt arXiv:math/0409010v1 [math.OC]}, 2004.

\bibitem{moreau05}
L.~Moreau.
\newblock Stability of multi-agent systems with time-dependent communication
  links.
\newblock {\em IEEE Transactions on Automatic Control}, 50:169--182, 2005.

\bibitem{olfati07}
R.~Olfati-Saber, J.A. Fax, and R.M. Murray.
\newblock Consensus and cooperation in networked multi-agent systems.
\newblock {\em Proceedings of the IEEE}, 95:215--233, 2007.

\bibitem{pecora98}
L.M. Pecora and T.L. Carroll.
\newblock Master stability functions for synchronized coupled systems.
\newblock {\em Physical Review Letters}, 80:2109--2112, 1998.

\bibitem{pogromsky01}
A.~Pogromsky and H.~Nijmeijer.
\newblock Cooperative oscillatory behavior of mutually coupled dynamical
  systems.
\newblock {\em IEEE Transactions on Circuits and Systems-I}, 48:152--162, 2001.

\bibitem{ren05}
W.~Ren and R.W. Beard.
\newblock Consensus seeking in multiagent systems under dynamically changing
  interaction topologies.
\newblock {\em IEEE Transactions on Automatic Control}, 50:655--661, 2005.

\bibitem{sontag98}
E.~Sontag.
\newblock {\em Mathematical Control Theory: Deterministic Finite Dimensional
  Systems}.
\newblock Springer, 1998.

\bibitem{stan07}
G.-B. Stan and R.~Sepulchre.
\newblock Analysis of interconnected oscillators by dissipativity theory.
\newblock {\em IEEE Transactions on Automatic Control}, 52:256--270, 2007.

\bibitem{stoop00}
R.~Stoop, K.~Schindler, and L.A. Bunimovich.
\newblock Neocortical networks of pyramidal neurons: from local locking and
  chaos to macroscopic chaos and synchronization.
\newblock {\em Nonlinearity}, 13:1515--1529, 2000.

\bibitem{strogatz01}
S.H. Strogatz.
\newblock Exploring complex networks.
\newblock {\em Nature}, 410:268--276, 2001.

\bibitem{tsitsiklis86}
J.N. Tsitsiklis, D.P. Bertsekas, and M.~Athans.
\newblock Distributed asynchronous deterministic and stochastic gradient
  optimization algorithms.
\newblock {\em IEEE Transactions on Automatic Control}, 31:803--812, 1986.

\bibitem{aut7874}
S.E. Tuna.
\newblock Synchronizing continuous-time neutrally stable linear systems via
  partial-state coupling.
\newblock {\em \tt arXiv:0801.3185v1 [math.OC]}, 2008.

\bibitem{walker69}
T.J. Walker.
\newblock Acoustic synchrony: two mechanisms in the snowy tree cricket.
\newblock {\em Science}, 166:891--894, 1969.

\bibitem{wang02}
X.F. Wang.
\newblock Complex networks: topology, dynamics and synchronization.
\newblock {\em International Journal of Bifurcation and Chaos}, 12:885--916,
  2002.

\bibitem{wu05con}
C.W. Wu.
\newblock Algebraic connectivity of directed graphs.
\newblock {\em Linear and Multilinear Algebra}, 53:203--223, 2005.

\bibitem{wu05}
C.W. Wu.
\newblock Synchronization in networks of nonlinear dynamical systems coupled
  via a directed graph.
\newblock {\em Nonlinearity}, 18:1057--1064, 2005.

\bibitem{wu95}
C.W. Wu and L.O. Chua.
\newblock Synchronization in an array of linearly coupled dynamical systems.
\newblock {\em IEEE Transactions on Circuits and Systems-I}, 42:430–447, 1995.

\end{thebibliography}
\end{document}